\documentclass[a4paper,10pt]{article}
\usepackage{amsmath}
\usepackage{amssymb}
\usepackage{bbm}
\usepackage[pdftex]{graphicx}
\usepackage{lineno}
\title{New results on affine invariant points}
\author{Olaf Mordhorst \footnote{Supported by 'DAAD Doktorandenstipendium'}}

\begin{document}

\maketitle
\numberwithin{equation}{section}

\newtheorem{thm}{Theorem}[section]
\newtheorem{cor}[thm]{Corollary}
\newtheorem{lem}[thm]{Lemma}
\newtheorem{prop}[thm]{Proposition}
\newtheorem{ax}{Axiom}

\newtheorem{defn}{Definition}[section]
\newtheorem{rem}{Remark}[section]

\begin{abstract}
We prove a conjecture of B. Gr\"unbaum stating that the set of affine invariant points of a convex body equals to the set of points invariant under all affine linear symmetries of the convex body. As a consequence we give a short proof on the fact that the affine space of affine linear points is infinite dimensional. In particular, we show that the set of affine invariant points with no dual is of second category. We investigate extremal cases for a class of symmetry measures. We show that the center of the John and L\"owner ellipsoid can be far apart and we give the optimal order for the extremal distance of the two centers.  
\end{abstract}
\maketitle

\section{Introduction}
In his seminal paper on symmetries of convex bodies \cite{Gruenbaum}, B. Gr\"unbaum addressed some questions about the symmetry structure of convex bodies and affine invariant points. Recently, M. Meyer, C. Sch\"utt and E. Werner answered some of these questions and extended the theory of affine invariant points (cf. \cite{SymMeasures}, \cite{AIP} and \cite{DAIP}). In this paper, we give some solutions to problems posed in these articles. Let us start with the notion of affine invariant points. 

We refer by \(d\) to the dimension of the convex body and we always assume \(d\geq 2\). An affine invariant point is a map \(p:\mathcal{K}_d\rightarrow \mathbb{R}^d\) from the set of d-dimensional convex bodies to \(\mathbb{R}^d\) such that for every convex body \(C\) and every invertible affine map \(T\) we have
\[
p(T(C))=T(p(C))\quad .
\]
We also require \(p\) to be continuous with respect to the Hausdorff distance and the euclidean norm. Well-known examples of affine invariant points are the centroid \(g\) (i.e. center of gravity), the Santal\'{o} point \(s\), the center of the John ellipsoid \(j\) and the center of the L\"owner ellipsoid \(l\) (see \cite{Gruenbaum}).

The main result of this paper is a complete answer to a question of Gr\"unbaum whether we can relate the symmetry structure of a convex body to the set of affine invariant points. Denote by 
\[
\mathfrak{P}_d(C)=\{p(C): p \text{ affine invariant point}\}
\] 
and by 
\[
\mathfrak{F}_d(C)=\{x\in\mathbb{R}^d: Tx=x\text{ for every T affine linear with } T(C)=C\}
\]
then we have
\begin{thm}\label{GruenbaumConjecture}
For every \(C\in\mathcal{K}_d\), we have
\[
\mathfrak{P}_d(C)=\mathfrak{F}_d(C)\quad.
\]
\end{thm}
There have been some major contributions to prove this conjecture. 
In \cite{AIP}, there is a proof for the case that the affine dimension of \(\mathfrak{P}_d(C)\) equals to \(d-1\). Recently, I. Iurchenko showed that this conjecture is true for the case \(\mathfrak{F}_d(C)=\mathbb{R}^d\) (see \cite{Iurchenko}). 
In \cite{Kuchment} (see also \cite{KuchmentEnglish} for an English translation), P. Kuchment proved that the conjectures holds for classes of points \(p: \mathcal{K}_d\rightarrow\mathbb{R}^d\) with \(p(T(C))=Tp(C)\) for all \(T\) affine map such that the linear part is an element of a fixed compact subgroup of \(GL(d)\). He obtained an even stronger variant of this conjecture for points with \(p(C)\in\mathrm{int}(C)\) which we prove analogously for affine invariant point (see Proposition \ref{ProperGruenbaum}). He proved the conjecture also for similarity invariant points by reducing the problem to convex bodies with volume \(1\) which leads to the compact case. The proof of Theorem \ref{GruenbaumConjecture} is based on ideas of \cite{Kuchment}: We reduce the problem to convex bodies in John position which leads also to just consider the compact subgroup \(O(d)\subseteq GL(d)\). We decided to avoid the fibre bundle language of \cite{Kuchment} and to replace the short and elegant but abstract topological arguments by explicit computations in order to make the proof more transparent to the intended audience. 

The set \(\mathfrak{P}_d\) of all affine invariant points on \(\mathcal{K}_d\) is an affine subspace of \(C(\mathcal{K}_d,\mathbb{R}^d)\), the space of continuous functions from \(\mathcal{K}_d\) to \(\mathbb{R}^d\), since one can show that every affine combination of affine invariant points is itself an affine invariant point. The following theorem answers a question of Gr\"unbaum and is proven in \cite{AIP}.
\begin{thm}
The affine dimension of \(\mathfrak{P}_d\) is infinite.
\end{thm} 

With Theorem \ref{GruenbaumConjecture} we provide an alternative proof of this theorem in a more concise way without using the technique of floating bodies.

It is a classical fact of convex geometry that the centroid and the Santal\'{o} point as well as the L\"owner and John point have a duality property in the following sense. L.A. Santal\'{o} proved in \cite{SantaloPoint} that we always have 
the identity \(s((C-g(C))^{o})=0\) where \(K^{o}\) denotes the dual of \(K\). For the L\"owner and John point we have
\(j((C-l(C))^{o})=0\) as well. A general definition for duality following \cite{DAIP} is that an affine invariant point \(q\) is dual to an affine invariant point \(p\) iff
\[
q((C-p(C))^{o})=0
\]
for every \(C\in\mathcal{K}_d\) (and providing the left-hand side is well-defined).
\(q\) is dual to \(p\) iff \(p\) is dual to \(q\) and we refer the reader to \cite{DAIP} for further properties on duality. A central result of this paper is the fact that there are affine invariant points with no dual. With Theorem \ref{GruenbaumConjecture}, we are able to show something even stronger. We equip \(\mathfrak{P}_d\) with the following distance
\[
 \mathrm{dist}(p_1,p_2)=\sup\limits_{B_2^d\subseteq C \subseteq d B_2^d}||p_1(C)-p_2(C)||_2
\]
(see also \cite{AIP}, section 3.2) and we have then
\begin{thm}
There exists an open and dense set \(W\subseteq\mathfrak{P}_d\) of affine points with no dual.
\end{thm}

In \cite{SymMeasures}, the following class of symmetry measures is introduced. For  \(p_1, p_2\in\mathfrak{P}_d\) such that we have \(p_1(C), p_2(C)\in C\) for every convex body \(C\), put

\begin{align}
\phi_{p_1, p_2}(C)=
\begin{cases}
 1 &\text{ if }p_1(C)=p_2(C)\\
 1-\frac{\|p_1(C)-p_2(C)\|_2}{\mathrm{vol}_1(a\cap C)}
 &\text{ else}
\end{cases}\notag
\end{align}
where \(a\) is the line through \(p_1(C)\) and \(p_2(C)\). The idea is that a convex body lacks symmetry if \(\phi\) is close to \(0\). The concern of \cite{SymMeasures} is to give extremal cases for this class of symmetry measures. We answer the open question if \(\phi_{p_1,p_2}(C)=0\) can happen for some \(p_1, p_2\) and we improve
a result for the case \(p_1=j\) and \(p_2=l\).
\subsection{Organistaion of the paper}
In Chapter 2, we give some definitions and background information which are important for further understanding of this paper. In Chapter 3, we prove the conjecture of Gr\"unbaum that \(\mathfrak{P}_d(C)=\mathfrak{F}_d(C)\) since this result will be a major ingredient for the subsequent investigations. In chapter 4, we give a minor extension of Theorem \ref{GruenbaumConjecture} with which we can show that
\(\mathfrak{P}_d\) is infinite dimensional. As a second application of Theorem \ref{GruenbaumConjecture}, we prove that the set of affine invariant points with no dual is of second category in \(\mathfrak{P}_d\). In the last chapter, we look at some questions about symmetry measures, especially, that the John and L\"owner point can be far apart to each other.

\section{Preliminaries}
For background information on affine invariant points, we refer the reader to \cite{Gruenbaum}, \cite{Kuchment}, \cite{KuchmentEnglish}, \cite{SymMeasures}, \cite{AIP} and \cite{DAIP}. We follow the presentation of \cite{AIP}, \cite{DAIP}. To be self-contained, we recall some notions found in these papers.
A d-dimensional convex body is a subset of \(\mathbb{R}^d\) which is convex,
compact and with non-empty interior. We denote by \(\mathcal{K}_d\) the set of d-dimensional convex bodies. A classical example for symmetric
convex bodies are the \(l_p^d\)-balls which are defined by

\begin{align}
B_p^d=\{x\in \mathbb{R}^d: \sum\limits_{k=1}^d |x_k|^p\leq 1\}\notag
\end{align}
for \(p\in [1,\infty)\) and

\begin{align}
B_{\infty}^d=\{x\in\mathbb{R}^d: \max\limits_{k=1\dots d} |x_k|\leq 1\}\quad.\notag
\end{align}
We equip \(\mathcal{K}_d\) with the
Hausdorff distance \(d_H\) which is for \(C_1, C_2 \in \mathcal{K}_d\) defined by

\begin{align}
d_H(C_1, C_2)
=\min\{\alpha\geq 0 : C_1\subseteq C_2+\alpha B_2^d \text{ and } 
C_2\subseteq C_1 + \alpha B_2^d\} \quad.\notag
\end{align}
We denote by \(GL(d)\) the group of linear, invertible  operators from \(\mathbb{R}^d\) to
\(\mathbb{R}^d\) and by \(O(d)\subseteq GL(d)\) the subgroup of orthogonal operators.
Let

\begin{align}
\mathcal{AT}_d=\{T: T=L+b, L\in GL(d) \text{ and }b\in\mathbb{R}^d\}\notag
\end{align}
be the set of affine linear and invertible transformations of \(\mathbb{R}^d\) where \(L+b\) should be understood as the map \(L+b:\mathbb{R}^d\rightarrow\mathbb{R}^d\), \((L+b)(x)=L(x)+b\).
From now on, we omit the \(d\) if it is clear from the context.

\begin{defn}
Let \(p:\mathcal{K}_d\rightarrow \mathbb{R}^d\) be a map which is continuous with respect to the Hausdorff distance and the euclidean norm. We call \(p\) an affine invariant point if  
\begin{align}
p(T(C))=T(p(C))\notag
\end{align} 
holds for every \(C\in\mathcal{K}_d\) and \(T\in\mathcal{AT}\).
\end{defn} 
This notion was introduced by Gr\"unbaum in \cite{Gruenbaum} and we adapt it here. It should be noted that it would be more adequate to call those maps affine equivariant points and throughout this paper, we refer by equivariance to the identity \(p(T(C))=Tp(C)\). Also, we call affine invariant points more shortly affine points.

Let \(\mathfrak{P}_d\) be the set of affine points. This is an affine subspace of the
vector space \(C(\mathcal{K}_d, \mathbb{R}^d)\) of continuous functions from \(\mathcal{K}_d\)
to \(\mathbb{R}^d\).
We put \(\mathfrak{P}_d(C)=\{p(C): p\in\mathfrak{P}_d\}\) and 

\begin{align}
\mathfrak{F}_d(C)=\{x\in\mathbb{R}^d: T(x)=x\text{ for every }T\in \mathcal{AT} \text{ with }T(C)=C\}\notag
\end{align}  
which are both affine subspaces of \(\mathbb{R}^d\). This means \(\mathfrak{F}_d(C)\) is the affine subspace of all points which do not change under symmetries of \(C\), i.e. under \(T\in\mathcal{AT}\) with \(T(C)=C\). For \(p\in\mathfrak{P}_d\), \(C\in\mathcal{K}_d\) and
\(T\in\mathcal{AT}\) with \(T(C)=C\) we have

\begin{align}
p(C)=p(T(C))=T(p(C))\notag
\end{align}
and hence, \(p(C)\in\mathfrak{P}_d(C)\), i.e., we always have \(\mathfrak{P}_d(C)\subseteq\mathfrak{F}_d(C)\). 
\\
For \(M\subseteq \mathcal{AT}\) and \(X\subseteq \mathcal{K}_d\), we put

\begin{align}
M(X)=\{T(C): T\in M\text{ and }C\in X\}\quad.\notag
\end{align}
For \(t\in\mathbb{R}\) and \(C\in\mathcal{K}_d\), we put
\((t,C)=\{(t,x): x\in C\}\).
\\
It is a classical fact of convex geometry that for every convex body \(C\) there exists a unique ellipsoid \(\mathcal{J}(C)\) of maximal volume inside C, called the John ellipsoid, and a unique ellipsoid \(\mathcal{L}(C)\) of minimal volume including \(C\), called the L\"owner ellipsoid. The maps \(\mathcal{J}:\mathcal{K}_d\rightarrow\mathcal{K}_d\) and
\(\mathcal{L}:\mathcal{K}_d\rightarrow\mathcal{K}_d\) are continuous with respect to the Hausdorff distance and we have \(\mathcal{J}(T(C))=T(\mathcal{J}(C))\) and \(\mathcal{L}(T(C))=T(\mathcal{L}(C))\) for every \(T\in\mathcal{AT}\) and \(C\in\mathcal{K}_d\). 
We say that a convex body is in John position (L\"owner position) if the euclidean ball is the John ellipsoid (L\"owner ellipsoid).
The following theorem is due to F. John and K. Ball (see \cite{Ball}, \cite{John}).
\begin{thm}\label{ContactPoints}
Let \(C\) be a convex body with \(B_2^d\subseteq C\). Then \(C\) is in John position iff there exist contact points \(v_1, \dots, v_m\in\partial C\cap\partial B_2^d\) and weights \(c_1, \dots, c_m>0\) with
\(
\sum\limits_{i=1}^m c_i v_i v_i^{tr}=\mathrm{Id}_d
\)
and
\(
\sum\limits_{i=1}^m c_iv_i=0_d
\).
\end{thm}
An analogue result is also true for the L\"owner position.
The following corollary is found in \cite{John}.
\begin{cor}\label{JohnCorollary}
Let \(C\in\mathcal{K}_d\) be a convex body in John position. Then we have \(C\subseteq dB_2^d\).
\end{cor}

\section{Proof of Gr\"unbaum's conjecture}
 Let
 \(\mathcal{K}_d^J=\{K\in \mathcal{K}_d: B_2^d\subseteq K\text{ John-ellipsoid of }K\}\)
 be the set of convex bodies in John position. \(\mathcal{K}_d^J\) is an \(O(d)\)-invariant subset of \(\mathcal{K}_d\),
 i.e. \(O(d)(\mathcal{K}_d^J)\subseteq \mathcal{K}_d^J\) 
 (and obviously, we may replace '\(\subseteq\)' by '\(=\)'). Moreover, \(\mathcal{K}_d^J\) is a closed subset of
 \(\mathcal{K}_d\). Let 
 \(\tilde{p}: \mathcal{K}_d^J\rightarrow \mathbb{R}^d\notag\)
be a continuous map with \(\tilde{p}(L(C))=L\tilde{p}(C)\) for every \(C\in\mathcal{K}_d^J\) and \(L\in O(d)\). Then we may extend \(\tilde{p}\) (in a unique way) to an affine point. We remark that for every \(C\in\mathcal{K}_d\) there are \(T\in\mathcal{AT}\)
and \(K\in\mathcal{K}_d^J\) with \(C=T(K)\).

\begin{lem}\label{ExtensionLemma}
  The map \(p:\mathcal{K}_d\rightarrow \mathbb{R}^d\), \(p(TK)=T\tilde{p}(K)\) for \(K\in\mathcal{K}_d^J\) and \(T\in \mathcal{AT}\)
  is an affine point.
 \end{lem}
 
\textit{Proof.} First, we show that p is well-defined. Let \(K, K'\in \mathcal{K}_d^J\)
 be two convex bodies in John position and let \(T, T'\) be two invertible affine maps such that
 \(T(K)=T'(K')\). We show \(T\tilde{p}(K)=T'\tilde{p}(K')\). From \(T'^{-1}T(K)=K'\), we deduce 
 \(B_2^d=J(K')=\mathcal{J}(T'^{-1}T(K))=T'^{-1}T\mathcal{J}(K)=T'^{-1}TB_2^d\) and hence, \(T'^{-1}T\in O(d)\). It follows

\begin{align}
  T'\tilde{p}(K')=T'\tilde{p}((T'^{-1}T)(K))=T'(T'^{-1}T)\tilde{p}(K)=T\tilde{p}(K)\notag
\end{align}
Next, we show that p is continuous. We show this by proving that for every convergent sequence
 \(\lim\limits_{n\rightarrow\infty} C_n=C\), there is a subsequence \((C_{n_m})_m\) such that 
 \(\lim\limits_{m\rightarrow\infty}p(C_{n_m})=p(C)\).
Let \(C_n=L_n(K_n)+b_n\) with \(K_n\in\mathcal{K}_d^J\), \(L_n\in GL(d)\) and \(b_n\in\mathbb{R}^d\)
 for \(n\in\mathbb{N}\) and \(C=L(K)+b\) with \(L, K, b\), accordingly. Since the John point \(j(C_n)=b_n\)
 is continuous as an affine point we have \(\lim\limits_{n\rightarrow\infty}b_n = b\) and hence, 
 \(\lim\limits_{n\rightarrow\infty} L_n(K_n)=L(K)\). There are constants \(\alpha_1, \alpha_2>0\) such that
 \(\alpha_1B_2^d\subseteq L(K)\subseteq \alpha_2 B_2^d\). Let \(0<\varepsilon<\alpha_1\) be fixed.
  Then there is an \(N\in \mathbb{N}\) such that for all \(n\geq N\), we have

\begin{align}
   (\alpha_1-\varepsilon)B_2^d\subseteq L_n(K_n)\subseteq (\alpha_2+\varepsilon)B_2^d\quad.\notag
\end{align}
We have \(B_2^d\subseteq M \subseteq dB_2^d\) for \(M\in\mathcal{K}_d^J\) and we
 conclude that for \(n\geq N\), we have

\begin{align}
 (\alpha_1-\varepsilon)B_2^d\subseteq dL_n(B_2^d)\text{ and }
 L_n(B_2^d)\subseteq (\alpha_2+\varepsilon)B_2^d\quad,\notag
\end{align}
i.e., the singular values of \(L_n\) are bounded from below by
 \(\frac{\alpha_1-\varepsilon}{d}\) and from above by \(\alpha_2+\varepsilon\). We consider the set 
of all operators with singular values bounded from below by \(\frac{\alpha_1-\varepsilon}{d}\) and
from above by \(\alpha_2+\varepsilon\):

\begin{align}
\{UDV\in GL(d): U,V\in O(d) \text{ and } D=\mathrm{diag}(s_1,\dots, s_d), \frac{\alpha_1-\varepsilon}{d}\leq s_i \leq \alpha_2+\varepsilon\}\notag
\end{align} 
where \(\mathrm{diag}(a_1, \dots, a_d)\) denotes the diagonal matrix with diagonal entries \(a_1, \dots, a_d\). This set is compact as the image of the compact set

\begin{align}
O(d)\times\{\mathrm{diag}(s_1,\dots,s_2)\in GL(d): \frac{\alpha_1-\varepsilon}{d}\leq s_i \leq \alpha_2+\varepsilon\}\times O(d)\notag
\end{align}
with respect to the continuous map \(\omega:GL(d)\times GL(d)\times GL(d)\rightarrow GL(d)\), \(\omega(A,B,C)=ABC\).
Therefore, there exists a convergent subsequence \((L_{n_m})_m\) with limit \(L'\) in \(GL(d)\). From 
 \(\lim\limits_{m\rightarrow\infty}L_{n_m}(K_{n_m})=L(K)\) it follows 
 \(\lim\limits_{m\rightarrow\infty}K_{n_m}=L'^{-1}L(K)\) and since \(\mathcal{K}_d^J\) is closed, we have
 \(L'^{-1}L(K)\in \mathcal{K}_d^J\) and hence, \(L'^{-1}L\in O(d)\). This yields:

 \begin{align}
  &\lim\limits_{m\rightarrow\infty}p(L_{n_m}(K_{n_m}))
 =\lim\limits_{m\rightarrow\infty}L_{n_m}\tilde{p}(K_{n_m})=L'\tilde{p}((L'^{-1}L)(K))\notag\\
 =& L'(L'^{-1}L)\tilde{p}(K)=L\tilde{p}(K)=p(L(K))\notag
 \end{align}
The proof of \(Tp(C)=p(T(C))\) is straight-forward.

\hfill \(\Box\)
\\[2ex]
We can now prove Gr\"unbaum's conjecture.
We only have to show that for
every \(x_0\in \mathfrak{F}_d(K)\) there is an affine point \(p\) such that \(p(K)=x_0\) since we always have \(\mathfrak{P}_d(K)\subseteq\mathfrak{F}_d(K)\). 

\begin{thm}\label{GruenbaumTheorem}
 Let \(K\in \mathcal{K}_d\) be a convex body, \(x_0\in \mathfrak{F}_d(K)\), then there is an affine point
 \(p: \mathcal{K}_d:\rightarrow \mathbb{R}^d\) such that \(p(K)=x_0\).
\end{thm}
\textit{Proof.} By an affine transformation of \(K\) and \(x_0\) we may assume without loss of generality that
\(K\) is in John position. Using the preceding lemma it is sufficient to construct a continuous and \(O(d)\)-equivariant
map \(\tilde{p}:\mathcal{K}_d^J\rightarrow \mathbb{R}^d\) such that \(\tilde{p}(K)=x_0\).
We construct \(\tilde{p}\) by an averaging argument. We start with the map
\(\theta: O(d)(K)\rightarrow \mathbb{R}^d\)       

\begin{align} 
\theta\left(L(K)\right)= L(x_0)\qquad.\notag
\end{align}
This map is well-defined and continuous: Suppose \(L(K)=L'(K)\), i.e. 
\(L'^{-1}L(K)=K\).
Since \(x_0\in \mathfrak{F}_d(K)\) it follows \(L'^{-1}L(x_0)=x_0\), i.e. \(L(x_0)=L'(x_0)\) and
hence \(\theta\) is well-defined. 

For the continuity we show that every convergent sequence in \(O(d)(K)\) has a subsequence such that
the image of this subsequence converges to the image of the limit. Suppose \(\left(L_n(K)\right)_{n\in\mathbb{N}}\) 
converges to an \(L(K)\in O(d)(K)\). By similar arguments as the ones in Lemma \ref{ExtensionLemma}, we may extract a 
convergent subsequence \(\left(L_{n_m}\right)_{m\in\mathbb{N}}\) with limit say \(L'\).
First, note that \(\lim_{m\rightarrow\infty}L_{n_m}(x_0)=L'(x_0)\).
From

\begin{align}
 L(K)=\lim\limits_{m\rightarrow\infty}L_{n_m}(K)=L'(K)\notag
\end{align}
we conclude that \(L'^{-1}L(K)=K\) and therefore, \(L'(x_0)=L(x_0)=\theta(L(K))\).
We remark that \(O(d)(K)\) is compact since \(\iota: O(d)\rightarrow\mathcal{K}_d\), 
\(\iota(L)=L(K)\) is continuous and therefore, \(O(d)(K)\) is compact as an image of a compact set, and in particular,
closed in \(\mathcal{K}_d^J\). By the Tietze-Urysohn lemma, there is a continuous extension of \(\theta\)
to a function \(\theta'\) on \(\mathcal{K}_d^J\).
Let \(\mu\) be the normalized (i.e. probability) Haar measure on \(O(d)\). We put \(\tilde{p}:\mathcal{K}_d^J\rightarrow\mathbb{R}^d\)

\begin{align}
 \tilde{p}(K)=\int\limits_{O(d)}L^{-1}\theta'(L(K))\mathrm{d}\mu(L)\notag
\end{align}
and show that \(\tilde{p}\) is continuous, \(O(d)\)-equivariant and \(\tilde{p}(K)=x_0\).
Indeed,

\begin{align}
\tilde{p}(K)=\int\limits_{O(d)}L^{-1}\theta'(L(K))\mathrm{d}\mu(L)=\int\limits_{O(d)}L^{-1}L(x_0)\mathrm{d}\mu(L)
\int\limits_{O(d)}x_0\mathrm{d}\mu(L)=x_0\quad.\notag
\end{align}
For every \(\Lambda\in O(d)\) and \(C\in\mathcal{K}_d^J\) we have using the invariance of the Haar measure (for \(O(d)\), left-
and right-invariance are the same):

\begin{align}
 \tilde{p}\left(\Lambda(C)\right)=&\int\limits_{O(d)}L^{-1}\theta'(L\Lambda(C))\mathrm{d}\mu(L)
=\int\limits_{O(d)}(L\Lambda^{-1})^{-1}\theta'((L\Lambda^{-1})k(C))\mathrm{d}\mu(L)\notag\\
=&\Lambda\int\limits_{O(d)}L^{-1}\theta'((L(C)))\mathrm{d}\mu(L)=\Lambda\tilde{p}(C)\notag
\end{align}
As for the continuity suppose the sequence \(\left(C_n\right)_{n\in\mathbb{N}}\) converges
to \(C\) in \(\mathcal{K}_d^J\). We have pointwise convergence of the integrand, i.e. \(L^{-1}\theta'(L(C_n))\)
converges to \(L^{-1}\theta'(L(C))\) for every \(L\in O(d)\). Furthermore, by compactness of \(\mathcal{K}_d^J\),
we have 

\begin{align}
B=\sup\limits_{M\in\mathcal{K}_d^J}||\theta'(M)||_2<\infty\notag
\end{align}
and we conclude 
\(\lim\limits_{n\rightarrow\infty}\tilde{p}(C_n)=\tilde{p}(C)\) using the theorem of dominated convergence.

\hfill \(\Box\) 

\section{The affine dimension of $\mathfrak{P}_d$ is infinite}
We show a slight extension of Theorem \ref{GruenbaumTheorem} in order to prove infinite dimensionality. We call a convex body \(C_1\) affine equivalent to a convex body \(C_2\) if there is a \(T\in\mathcal{AT}\)
with \(T(C_1)=C_2\). One checks immediately, that this is an equivalence relation and from now on, we just say \(C_1\) and \(C_2\) are affinely equivalent.

\begin{cor}\label{ExtendedGruenbaum}
 Let \(K_1, \dots, K_m\) be convex bodies which are pairwise not affinely equivalent, and let \(x_i\in \mathfrak{F}(K_i)\) for \(1\leq i \leq m\). Then there is an
affine point \(p\) such that \(p(K_i)=x_i\) for \(1\leq i \leq m\).
\end{cor}
\textit{Proof.} The proof relies on similar techniques as the proof of Theorem \ref{GruenbaumTheorem}. Assume without loss of generality that all convex bodies concerned are in John position. Put

\begin{align}
\theta: \bigcup_{i=1}^m O(d)(K_i)&\rightarrow \mathbb{R}^d\notag\\
L(K_j)&\mapsto L(x_j)\quad.\notag
\end{align}
The sets \(O(d)(K_i)\) are pairwise disjoint, since no two of the \(K_i\)'s are affinely equivalent. By compactness of the \(O(d)(K_i)\), these sets have positive distance with respect to the Hausdorff distance and hence, we can conclude from the continuity of \(\theta|_{O(d)(K_i)}\)
that \(\theta\) is itself continuous. Using this \(\theta\), we now proceed as in the proof of Theorem \ref{GruenbaumTheorem}.
\hfill \(\Box\) 
\\[2ex]
Using this corollary, we can construct a sequence of affine points which are linearly independent and hence, also affinely independent.

\begin{thm}
The affine dimension of\(\mathfrak{P}_d\) is infinite.
\end{thm}
\textit{Proof.} Assume, we have a sequence of convex bodies \(K_k\), \(k=1,2,3, \dots\) which are pairwise not affinely equivalent and such that \(\mathfrak{F}_d(K_k)\) does not reduce to one point for every \(k\). Start with an arbitrary affine point \(p_1\). The sequence only containing \(p_1\) is of course linearly independent. Now, assume we have a sequence \(p_1, \dots, p_n\) of linearly independent affine points. For every choice of \(x_k\in \mathfrak{F}_d(K_k)\), \(k=1, \dots, n+1\), we find by Theorem \ref{ExtendedGruenbaum} an affine point \(p\) with \(p(K_k)=x_k\) for \(k=1,\dots, n+1\). In particular, we find an affine point \(p_{n+1}\) linearly independent to the sequence \(p_1,\dots, p_n\). Hence, we obtain a sequence of infinite many linearly independent affine points and therefore, \(\mathfrak{P}_d\) is infinite dimensional. 
The only problem left open is the question if such sequences of \(K_k\)'s exist. This should be clear at least by the fact, that polytopes with different numbers of vertices cannot be affinely equivalent because affine, invertible operators always map vertices to vertices. 
However, we work out an explicit construction at least for the case \(d\geq 3\).
For \(k=1,2,3, \dots\), put

\begin{align}
K_k=\mathrm{conv}\left[\{0_{d}\}, (1, B_{k+1}^{d-1})\right]\notag
\end{align}
where \(0_{d}\in\mathbb{R}^{d}\) is the origin. 
The \(K_k\)'s are cones with base \(B_{k+1}^{d-1}\) and the construction of a double cone over a \((d-1)\)-dimensional convex body was already considered in [\cite{Kuchment}, Theorem 5] for a somewhat related problem on similarity invariant points. For no two \(k, k'\) with \(k\neq k'\), the cones \(K_k\) and \(K_{k'}\) are affine equivalent.
The set of extreme points of \(K_k\) is 
\(\{0_{d}\}\cup\{0\}\times\partial B_{k+1}^{d-1}\) (note, that \(k+1>1\)). Since every affine invertible operator maps extreme points to extreme points, we conclude that every affine operator \(T\) with \(T(K_k)=K_k\) maps \(0_d\) to \(0_d\) and \(\{1\}\times\partial B_{k+1}^{d-1}\) to \(\{1\}\times\partial B_{k+1}^{d-1}\) and hence, \(T\) does not change the first coordinate. 
Hence, we deduce that \(\mathfrak{F}_{d}(K_k)\supseteq \mathbb{R}\times\{0_{d-1}\}\). \(K_k\) is invariant under permutations of the last \(d-1\) coordinates and that yields \(\mathfrak{F}_{d+1}(K_k)=\mathbb{R}\times \{0_{d-1}\}\). 

\hfill \(\Box\)

\section{Dual affine invariant points}
We call an affine point \(p\) proper if for every \(C\in\mathcal{K}_d\), we have \(p(K)\in \mathrm{int}(K)\). Let \(p,q\) be two proper affine points. We say that \(q\) is dual to \(p\) if \(q((C-p(C))^{0})=0\) for every \(C\in\mathcal{K}_d\). In \cite{DAIP}, it was pointed out that there are proper affine points with no dual. Here, we want to show that the set of affine points with no dual is of second category with respect to some distance. More precisely, we show that there is an open and dense set of affine points with no dual.
We make use of the following fact
(cf. \cite{DAIP}, Lemma 3.6, Theorem 4.3):

\begin{thm}
 A proper affine point has a dual iff for every \(C\in\mathcal{K}_d\), there is at most
 one \(z\in int(C)\) with \(p((C-z)^{\circ})=0\).
\end{thm}
From [\cite{DAIP}, Proposition 3.13] we know that there is always at least one \(z\) with \(p((C-z)^{\circ})=0\).

We use the following metric to measure the distance of two affine points which was introduced
in [\cite{AIP}, section 3.2]:

\begin{align}
 \mathrm{dist}(p,q)=\sup\limits_{B_2^d\subseteq C \subseteq d B_2^d}||p(C)-q(C)||_2\notag
\end{align}
With respect to this distance we have the following topological result on proper affine points with no dual.
\begin{thm}
 There exists an open and dense set \(W\subseteq \mathfrak{P}_d\) of affine points with no dual.
\end{thm}
\textit{Proof.} Let \(p\) be an affine point having a dual and \(\varepsilon>0\). We construct an
affine point \(q\) such that \(\mathrm{dist}(p,q)<\varepsilon\) and such that there is a neighbourhood \(V\) of \(q\) of affine points with no dual.
We will construct \(V\) in such a way that for every \(q'\in V\) there is a
\(z_0\in \mathrm{int}(B_1^d)\backslash\{0\}\) with \(q'((B_1^d-z_0)^{\circ})=0\). Since we always have \(q(B_{\infty}^d)=0_d\)
we conclude from the preceding discussion that \(q'\) has no dual. 

Since \(p\) is continuous there exists a \(\eta>0\) such that for every \(C\in\mathcal{K}_d\)
with \(d_H(C, B_{\infty}^d)< \eta\) we have \(\|p(C)\|_2<\frac{\varepsilon}{4d}\). For every \(\delta\in(-1,1)\)
we can choose \(l_{\delta}\), \(c_{\delta}\) and \(a_{\delta}\) such that \(T_{\delta}((B_1^d-\delta e_1)^{\circ})\) is
in John position where \(T_{\delta}\) is the affine map

\begin{align}
T_{\delta}= 
\begin{pmatrix}
l_{\delta} & 0 & \dots & 0\\
0\\
\vdots & &c_{\delta}\cdot I_{d-1}\\
0
\end{pmatrix}
+
\begin{pmatrix}
a_{\delta}\\
0\\
\vdots\\
0\end{pmatrix}\quad.\notag
\end{align}
We have \(\lim_{\delta\rightarrow 0} l_{\delta}=1\), \(\lim_{\delta \rightarrow 0} c_{\delta}=1\) and 
\(\lim_{\delta\rightarrow 0} a_{\delta}=0\), moreover

\begin{align}
\lim\limits_{\delta\rightarrow 0}T_{\delta}\left((B_1^d-\delta e_1)^{\circ}\right)=B_{\infty}^d\quad.\notag
\end{align}
Therefore, we may choose sufficiently small \(0<\delta_1< \delta_2<1\) such that
for \(i=1,2\) we have
\begin{align}
\left|\frac{(-1)^i l_{\delta_i}\cdot\varepsilon}{5d}+a_\delta\right|< \frac{\varepsilon}{4d}\notag
\end{align}
and
\begin{align}
d_H(B_{\infty}^d, T_{\delta_i}\left((B_1^d-\delta_i e_1)^{\circ}\right))<\eta\quad.\notag
\end{align}
We put \(C_i=T_{\delta_i}\left((B_1^d-\delta_i e_1)^{\circ}\right)\).
We have \(\mathfrak{F}_d((B_1^d-\delta_i e_1)^{\circ})=\mathbb{R}\times\{0_{d-1}\}\). Note that \(B_{\infty}^d, (B_1^d-\delta_1 e_1)^{\circ} \text{ and } (B_1^d-\delta_2 e_1)^{\circ}\) are pairwise not affine equivalent, i.e. there exists no affine map which transforms one convex body into another.
Hence, we can choose an affine point \(r\)
such that 
\begin{align}
r(B_{\infty}^d)=0_d, \quad r((B_1^d-\delta_1 e_1)^{\circ})=-\frac{\varepsilon}{5d}e_1\text{ and } r((B_1^d-\delta_2 e_1)^{\circ})=\frac{\varepsilon}{5d}e_1\quad.\notag
\end{align} 
It follows that 
\begin{align}
\|r(C_i)\|_2=\left|\frac{(-1)^i l_{\delta_i}\cdot \varepsilon}{5d}+a_{\delta_i}\right|<\frac{\varepsilon}{4d}\quad.\notag
\end{align}
Since the compact sets \(\{B_{\infty}^d\}, O(d)(C_1)\text{ and }O(d)(C_2)\) are pairwise disjoint we conclude\begin{align}
\alpha_i=\inf\limits_{L\in O(d)}d_H(L(C_i), B_{\infty}^d)>0\quad.\notag
\end{align}
We put \(\alpha=\min[\alpha_1,\alpha_2]\).
 Let \(U_i\subseteq \mathcal{K}_d^J\) be relatively 
open neighbourhoods of \(C_i\) such that

\begin{align}
\|r(K)\|_2<\frac{\varepsilon}{4d},\quad
d_H(K, B_{\infty}^d)< \eta \text{ and }
d_H(K, C_i)<\frac{\alpha}{2}\notag
\end{align}
for every \(K\in U_i\) and put \(U=U_1\cup U_2\). We conclude that for every \(K\in U_i\) and \(L\in O(d)\) we have by the reverse triangle inequality

\begin{align}
d_H(B_{\infty}^d, L(K))&\geq d_H(B_{\infty}^d, L(C_i))-d_H(L(C_i), L(K))\notag\\&=d_H(B_{\infty}^d, L(C_i))-d_H(C_i, K)
> \alpha -\frac{\alpha}{2}=\frac{\alpha}{2}\quad.\notag
\end{align}
By the Tietze-Urysohn lemma there is a
continuous map \(\phi:\mathcal{K}_d^J\rightarrow [0,1]\)
with \(\phi(L(C_i))=1\) for \(L\in O(d)\) and \(\phi(K)=0\)
for every \(K\in (O(d)(U))^{c}\). Put 
\(\Phi:\mathcal{K}_d^J\rightarrow [0,1]\)

\begin{align}
\Phi(K)=\int\limits_{O(d)} \phi(L(K))\mathrm{d}\mu(L)\notag
\end{align}
where \(\mu\) is the normalized Haar measure on \(O(d)\).
Put \(\tilde{q}:\mathcal{K}_d^J\rightarrow \mathbb{R}^d\),
\(\tilde{q}(K)=\Phi(K)r(K)+(1-\Phi(K))p(K)\) which is continuous as a combination of continuous maps and \(O(d)\)-equivariant. Denote by \(q\) the unique extension of \(\tilde{q}\) to an affine point. For \(K\in \mathcal{K}_d^J\backslash O(d)(U)\) we have \(\|q(K)-p(K)\|_2=0\) and
for \(K\in U\) and \(L\in O(d)\) we have

\begin{align}
\|q(L(K))-p(L(K))\|_2&= \Phi(L(K))\|r(L(K))-p(L(K))\|_2
\leq \|L(r(K)-p(K))\|_2\notag\\
&=\|r(K)-p(K)\|_2\leq \|r(K)\|_2+\|p(K)\|_2\notag\\
&<\frac{\varepsilon}{4d}+\frac{\varepsilon}{4d}=\frac{\varepsilon}{2d}\quad.\notag
\end{align}
  
Let \(C\in\mathcal{K}_d\) be an arbitrary convex body with the property
\(B_2^d\subseteq C \subseteq dB_2^d\). Then there is an affine transformation \(T=L+b\), \(L\in GL(d)\), \(b\in\mathbb{R}\) and a convex body \(K\in\mathcal{K}_d^J\) with
\(C=T(K)\). We have the following estimate:

\begin{align}
\|q(C)-p(C)\|_2=&\|q(T(K))-p(T(K))\|_2=\|L(q(K)-p(K))\|_2\notag\\
\leq & \|L\|_{Op}\|q(K)-p(K)\|_2<\|L\|_{Op}\frac{\varepsilon}{2d}\notag
\end{align}
where \(\|\cdot\|_{Op}\) denotes the spectral norm.
We have \(b=b+L(0_d)\in b+L(K)\subseteq dB_2^d\) and hence,
\(dB_2^d-b\subseteq 2dB_2^d\), and we conclude

\begin{align}
L(B_2^d)\subseteq L(K)\subseteq dB_2^d-b\subseteq 2dB_2^d\notag
\end{align}
which leads to the bound of the spectral norm of \(L\) by \(2d\)
and therefore,

\begin{align}
\|q(C)-p(C)\|_2< \|L\|_{Op}\frac{\varepsilon}{2d}\leq \varepsilon\notag
\end{align}
or \(\mathrm{dist}(q,p)<\varepsilon\), respectively.
Put 
\begin{align}
V=\{q'\in\mathfrak{P}_d: q'((B_1^d-\delta_1 e_1)^{\circ})_1<0\text{ and }q'((B_1^d-\delta_2 e_1)^{\circ})_1>0\}\notag
\end{align}
where \(q'(\dots)_1\) is the first component of \(q'(\dots)\). The set \(V\) is an open neighbourhood of \(q\). Note that for every \(q'\in V\) the map
\begin{align}
\psi: [\delta_1,\delta_2]\rightarrow \mathbb{R}\times \{0_{d-1}\},\quad \psi(\delta)=q'\left((B_1^d-\delta e_1)^{\circ}\right)\notag
\end{align}
is continuous and hence, for every \(q'\) there is a \(\delta'\in[\delta_1,\delta_2]\) with
\begin{align}
q'\left((B_1^d-\delta'e_1)^{\circ}\right)=0e_1=0_d\notag
\end{align}
by the intermediate value theorem.
\hfill \(\Box\)

\section{Symmetry measures and the maximal distance of the John and L\"owner point}
When we talk about symmetry measures, we mean the class introduced in \cite{SymMeasures}. 
Consider two affine points \(p_1, p_2\) with \(p_1(C), p_2(C)\in C\) for every \(C\in \mathcal{K}_d\). We put
\(\delta_{p_1,p_2}:\mathcal{K}_d\rightarrow [0,1]\)

\begin{align}
\delta_{p_1, p_2}(C)=0 \text{, if }p_1(C)=p_2(C)\notag\\
\delta_{p_1, p_2}(C)=\frac{\|p_1(C)-p_2(C)\|}{\mathrm{vol}_1(a\cap C)}\notag
\end{align}
where \(a\) denotes the line through \(p_1(C), p_2(C)\) and define the symmetry measure as \(\phi_{p_1,p_2}=1-\delta_{p_1,p_2}\). Morally, if the value of this symmetry measure is close to \(1\), we have high symmetry. For example this value is always \(1\) for the simplex and for symmetric convex bodies. In \cite{SymMeasures}, the question is raised if \(\phi_{p_1,p_2}(C)=0\) can occur for some \(p_1, p_2, C\). In part 1, we answer this question in the affirmative for every dimension \(d\geq 2\).
In part two, we treat the special case where \(p_1=j\) and \(p_2=l\), i.e. the case of the John and L\"owner point.

\subsection{Extremal $p_1, p_2$ for $\phi$}
We start with two propositions how to construct proper affine points and more generally, affine points with \(p(C)\in C\) for every \(C\in\mathcal{K}_d\). The proof of the first proposition is almost analogous to the proof of the same statement about different classes of equivariant points considered in \cite{Kuchment} but we repeat the arguments in our setting for the sake of completeness.

\begin{prop}\label{ProperGruenbaum}
 Let \(K\in\mathcal{K}_d\) be a convex body and \(x_0\in\mathfrak{F}_d(K)\cap\mathrm{int}(K)\).
Then there is a proper affine point \(p\) with \(p(K)=x_0\)
\end{prop}
\textit{Proof.} Without loss of generality we assume \(K\) to be in the John position. By Theorem \ref{GruenbaumTheorem},
we have an \(O(d)\)-equivariant, continuous map \(\tilde{p}:\mathcal{K}_d^J\rightarrow \mathbb{R}^d\) with 
\(\tilde{p}(K)=x_0\). By continuity, we have an open neighbourhood \(U\) of \(K\) in \(\mathcal{K}_d^J\)
with \(\tilde{p}(C)\in\mathrm{int}(C)\) for every \(C\in U\). The same holds for the open \(O(d)\)-invariant set

\begin{align}
 W=\bigcup\limits_{l\in O(d)} l(U)\quad.\notag
\end{align}
By the Tietze-Urysohn lemma, there is a continuous function \(\phi:\mathcal{K}_d^J\rightarrow[0,1]\) such that

\begin{align}
 \phi(C)=
\begin{cases}
 1 & \text{ if }C\in O(d)(K)\\
 0 & \text{ if }C\in W^{c}
\end{cases}\quad.\notag
\end{align}
Let \(\mu\) be the normalized Haar measure on \(O(d)\) and define \(\Phi:\mathcal{K}_d^J\rightarrow [0,1]\) by

\begin{align}                   
                   \Phi(C)=\int\limits_{O(d)}\phi(L(K))\mathrm{d}\mu(L)\quad.\notag
\end{align}
Obviously, \(\Phi\) is a continuous function with \(1\) on \(O(d)(K)\) and \(0\) on \(W^{c}\). We put 
\(\tilde{q}:\mathcal{K}_d^J\rightarrow\mathbb{R}^d\)

\begin{align}
\tilde{q}(C)=\Phi(C)\tilde{p}(C)+(1-\Phi(C))g(C)\notag
\end{align}
where \(g\) is the centroid. Then, \(\tilde{q}\) is continuous, \(O(d)\)-equivariant and we have 
\(\tilde{q}(C)\in\mathrm{int}(C)\) for every \(C\in\mathcal{K}_d^J\). The same holds if we extend
\(\tilde{q}\) to an affine point \(q\).
\hfill\(\Box\)
\\[2ex]

\begin{prop}\label{InteriorPoints}
 Let \(K\) be a convex body and \(x_0\in\mathfrak{F}_d(K)\cap\partial K\). Then there exists an affine 
 point \(p\) with \(p(K)=x_0\) and \(p(C)\in C\) for every convex body \(C\).
\end{prop}
\textit{Proof.} Since \(\mathfrak{F}_d(K)\) is an affine subspace of \(\mathbb{R}^d\) and
\(\frac{1}{2}(g(K)+x_0)\in \mathrm{int}(K)\) we conclude from Proposition \ref{ProperGruenbaum} that there is a
proper affine point \(q\) with \(q(K)=\frac{1}{2}(g(K)+x_0)\). We put \(\gamma_q: \mathcal{K}_d\rightarrow [0,1]\)

\begin{align}
              \gamma_q(C)= \begin{cases}
                             0 \text{ if }q(C)=g(C)\notag\\
                             \left(\sup\{\lambda>0: g(C)+\lambda(q(C)-g(C))\in C\}\right)^{-1}\text{ else}\notag
                            \end{cases}
\end{align}
We omit the proof that \(\gamma_q\) does not change under affine transformations, i.e. \(\gamma_q(T(C))=\gamma_q(C)\) for affine 
maps \(T\), and \(\gamma_q\) is continuous. The function \(\gamma_q\) has the nice property that for
every \(C\) and \(\gamma\) with \(\gamma\geq \gamma_q(C)>0\) we have \(g(C)+\frac{1}{\gamma}(q(C)-g(C))\in C\) and especially,
\(g(K)+\frac{1}{\gamma_q(K)}(q(K)-g(K))\in\partial K\). Put \(\psi:[0,1]\rightarrow[\frac{1}{2},1]\)
with \(\psi(s)=1-s\) for \(s\leq \frac{1}{2}\) and \(\psi(s)=s\) for \(s\geq \frac{1}{2}\) and define 
\(\Gamma_q=\psi\circ \gamma_q\). \(\Gamma_q\) is itself continuous as a composition of continuous functions and does not change under affine transformations. 
We can then define the desired affine invariant point \(p\) as \(p(C)=g(C)+\frac{1}{\Gamma_q(C)}(q(C)-g(C))\)
for \(C\in \mathcal{K}_d\). 

\hfill\(\Box\)

\begin{rem}
 Using the centroid \(g\) in the proof is not crucial. Actually, we could have taken any other proper affine
point instead.
\end{rem}

\begin{thm}
For every \(d\geq 2\) there exists a symmetry measure \(\phi_{p_1, p_2}\) and a convex body \(C\in\mathcal{K}_d\) with \(\phi_{p_1,p_2}(C)=0\).
\end{thm}
\textit{Proof:} Choose a convex body \(K\in\mathcal{K}_d\) such that the affine dimension of \(\mathfrak{P}_d(K)\) is at least \(1\). Those convex bodies do always exist because the set of convex bodies \(C\) with \(\mathfrak{P}_d(C)=\mathbb{R}^d\) is dense by [\cite{AIP}, Theorem 3]. If we want to avoid this implicit argument, one can choose as \(C\) for example the cone
\[
\mathrm{conv}\left[\{0_d\}, (1, B_2^{d-1})\right]
\]
for the case \(d\geq 3\). Since the affine point \(g(K)\) is an element of \(\mathrm{int}(K)\), the set \(\mathfrak{P}_d(K)\cap \partial K\) has at least two distinctive points, say \(x_1, x_2\). By Proposition \ref{InteriorPoints}, we find two affine points \(p_1, p_2\) with \(p_1(C),p_2(C)\in C\) for every \(C\in\mathcal{K}_d\) and \(p_1(K)=x_1\) and \(p_2(K)=x_2\).
This yields \(\phi_{p_1,p_2}(C)=0\).

\hfill \(\Box\)
\subsection{An extremal case for the L\"owner and John point}
What we will prove here is that 

\begin{align}
\frac{2}{d+1}\leq\inf\limits_{C\in\mathcal{K}_d}\phi_{j,l}(C)= \frac{2}{d+1}(1+o(1))\quad.\label{NewOrderOfMagnitude}
\end{align}
In \cite{SymMeasures}, we find the estimate

\begin{align}
\frac{2}{d+1}\leq\inf\limits_{C\in\mathcal{K}_d}\phi_{j,l}(C)\leq \frac{1}{2}\quad.\label{OldOrderOfMagnitude}
\end{align}
It turns out that, even though the right-hand side of (\ref{OldOrderOfMagnitude}) is not the right order of magnitude, we only have to modify the construction of \cite{SymMeasures} slightly to obtain the right order \(\frac{2}{d+1}\). We start with the left-hand side of (\ref{NewOrderOfMagnitude}). Actually, you can find a proof in [\cite{SymMeasures}, section 3, remark]. We recall it for two reasons. The first reason is that this proof is done for the Santal\'{o} point and the centroid (although the arguments work analogously for the John and L\"owner point). Secondly, there is a minor mistake in the proof. 
\begin{lem}
Let \(d\geq 2\). For every convex body \(C\) we have \(\phi_{l,j}(C)\geq \frac{2}{d+1}\).
\end{lem}
\textit{Proof.} Let \(\mathcal{E}\) be the John ellipsoid of \(C-j(C)\). By Corollary \ref{JohnCorollary}, we have \(C-j(C)\subseteq d\mathcal{E}\).
We conclude \(C-j(C)\subseteq d\mathcal{E}=-d\mathcal{E}\subseteq d(j(C)-C)\). Similarly, we have \(C-l(C)\subseteq d(l(C)-C)\). Let \(a\) be an arbitrary chord in \(C-j(C)\) passing through the origin. Then \(a\) is split by the origin into two line segments \(a_1, a_2\). From \(C-j(C)\subseteq d(j(C)-C)\), we deduce that \(a_1\subseteq d a_2\) and hence, 

\begin{align}
\frac{\mathrm{vol}_1(a_2)}{\mathrm{vol}_1(a_1)}\geq\frac{1}{d}\qquad\text{ or equivalently, }\quad\frac{\mathrm{vol}_1(a_2)}{\mathrm{vol}_1(a)}\geq \frac{2}{d+1}\quad.\notag
\end{align} 
The same holds for every chord through the L\"owner point. Assume \(l(C)\neq j(C)\) and let now \(a\) be the chord of \(C\) through \(j(C)\) and \(l(C)\). By the preceding discussion, the portion of the part of the chord which is not covered by the line segment between \(l(C), j(C)\) is at least \(\frac{1}{d+1}+\frac{1}{d+1}=\frac{2}{d+1}\) of the length of \(a\).

\hfill\(\Box\)
\\[2ex]
Now, we establish the right-hand side of \ref{NewOrderOfMagnitude}. We need the following lemma.

\begin{lem}\label{LoewnerLemma}
There is a sequence \((\varepsilon_d)_{d\in\mathbb{N}}\) with \(\varepsilon_d=\frac{1}{d}+O(\frac{1}{d^2})\) such that the L\"owner ellipsoid of 
\[
K_d=\mathrm{conv}\left[(-\varepsilon_d,\sqrt{1-\varepsilon_d^2}B_2^d), (\sqrt{1-\frac{1}{d}},\frac{1}{\sqrt{d}}B_2^d)\right]
\]
is \(B_2^{d+1}\).
\end{lem}
\textit{Proof.} For fixed \(0<\varepsilon<1\) and \(d\) consider

\begin{align}
K=\mathrm{conv}\left[(-\varepsilon,\sqrt{1-\varepsilon^2}B_2^d), (\sqrt{1-\frac{1}{d}},\frac{1}{\sqrt{d}}B_2^d)\right]\quad.\notag
\end{align}
For \(2\leq i \leq d+1\) put

\begin{align}
v_i^{\pm}=\pm \sqrt{1-\varepsilon^2}e_i-\varepsilon e_1
\text{ and } w_i^{\pm}=\pm \frac{1}{\sqrt{d}}e_i+\sqrt{1-\frac{1}{d}}e_1\notag
\end{align}
where \(e_i\) is the i-th unit coordinate vector.
We note that \(v_i^{\pm}, w_i^{\pm}\in \partial K\cap \partial B_2^{d+1}\) and

\begin{align}
V&=\sum\limits_{o\in\{\pm\}, i}v_i^{o}{v_i^{o}}^{\mathrm{tr}}=2
\begin{pmatrix}
d\varepsilon^2& & &\textbf{0}\\
 & 1-\varepsilon^2\\
 & & \ddots \\
\textbf{0} & & & 1-\varepsilon^2
\end{pmatrix}\notag\\
W&=\sum\limits_{o\in\{\pm\}, i}w_i^{o}{w_i^{o}}^{\mathrm{tr}}=2
\begin{pmatrix}
d(1-\frac{1}{d}) & & & \textbf{0}\\
 & \frac{1}{d}\\
 & & \ddots \\
 \textbf{0} & & & \frac{1}{d}
\end{pmatrix}\quad.\notag
\end{align}
If we find \(t_1, t_2>0\) such that \(t_1V+t_2W=\mathrm{Id}_{d+1}\) and 
\[
t_1\left(\sum\limits_{o\in\{\pm\}, i}v_i^{o}\right)
+t_2\left(\sum\limits_{o\in\{\pm\},i}w_i^{o}\right)=0\quad,
\]
we may conclude that \(B_2^{d+1}\) is the L\"owner ellipsoid of \(K\) (see Theorem \ref{ContactPoints}). These two equalities are equivalent to the following three equations:

\begin{align}
t_1(1-\varepsilon^2)+t_2 \frac{1}{d}=\frac{1}{2}, \quad
t_1 d\varepsilon^2+t_2 d(1-\frac{1}{d})&=\frac{1}{2}\text{ and }
t_1\varepsilon-t_2\sqrt{1-\frac{1}{d}}&=0\quad.\notag
\end{align}
From \(t_1(1-\varepsilon^2)+t_2 \frac{1}{d}=t_1 d\varepsilon^2+t_2 d(1-\frac{1}{d})\) we get

\begin{align}
\frac{t_1}{t_2}=\frac{d-1-\frac{1}{d}}{1-(d+1)\varepsilon^2}\notag
\end{align}
and the third equations is equivalent to

\begin{align}
\frac{t_1}{t_2}=\frac{\sqrt{1-\frac{1}{d}}}{\varepsilon}\quad.\notag
\end{align}
We obtain that

\begin{align}
\sqrt{1-\frac{1}{d}}\left(1-(d+1)\varepsilon^2\right)
=(d-1-\frac{1}{d})\varepsilon\notag
\end{align}
and the only positive solution of this quadratic equation in \(\varepsilon\) is

\begin{align}
\varepsilon_d&=\frac{\sqrt{1+4(d+1)\frac{1-\frac{1}{d}}{d^2-d-1}}-1}{2(d+1)\frac{\sqrt{1-\frac{1}{d}}}{d-1-\frac{1}{d}}} = \frac{\sqrt{1+\frac{4}{d}+O(\frac{1}{d^2})}-1}{2+O(\frac{1}{d})}\notag\\
&= \frac{1}{d}+O\left(\frac{1}{d^2}\right)\label{EpsilonFormula}
\end{align}
With \(\varepsilon=\varepsilon_d\) we verify that we have positive solutions \(t_1, t_2\). 

\hfill \(\Box\)
\\[2ex]
Consider the convex body

\begin{align}
C_d=\mathrm{conv}\left[(-\varepsilon_d, \sqrt{1-\varepsilon_d^2}\Delta_d), (\sqrt{1-\frac{1}{d}}, \frac{1}{\sqrt{d}}B_2^d)\right]\notag
\end{align}
where \(\Delta_d\) is the regular \(d\)-dimensional simplex inscribed in \(B_2^d\). 

\begin{thm}
Asymptotically in \(d\) we have the estimate \(\phi_{j,l}(C_d)\leq \frac{2}{d+1}(1+o(1))\).
\end{thm}
\textit{Proof.} The convex body \(C_d\) is similar to the construction in  [\cite{SymMeasures}, section 5]. The only difference is another scaling factor for \(\Delta_d\). We also use the same arguments as in [\cite{SymMeasures}, proof of Proposition 15] to show that the L\"owner ellipsoid of \(C_d\) is the same as the L\"owner ellipsoid of \(K_d\) as defined in Lemma \ref{LoewnerLemma} and hence, the L\"owner ellipsoid is \(B_2^{d+1}\). We can also argue in the same way as in [\cite{SymMeasures}, Proposition 15] that the John ellipsoid of \(C_d\) is included in 

\begin{align}
\mathrm{conv}\left[(-\varepsilon_d,\frac{\sqrt{1-\varepsilon_d^2}}{d} B_2^d),  (\sqrt{1-\frac{1}{d}}, \frac{1}{\sqrt{d}}B_2^d)\right]\notag
\end{align}
and hence, they share the same John ellipsoid.
This is also the John ellipsoid of

\begin{align}
\mathrm{conv}\left[(-\varepsilon_d,\sqrt{1-\varepsilon_d^2} \Delta_d),  (\sqrt{1-\frac{1}{d}}, \sqrt{d}\Delta_d)\right]\quad.\notag
\end{align}
This set is the frustum or truncation of a simplex. The exact expression for this simplex is

\begin{align}
\tilde{\Delta}_{d+1}=\mathrm{conv}\left[- \rho_d e_1, (\sqrt{1-\frac{1}{d}}, \sqrt{d}\Delta_d)\right]\notag
\end{align}
where we put

\begin{align}
\rho_d=\frac{\varepsilon_d\sqrt{d}+\sqrt{1-\varepsilon_d^2}\sqrt{1-\frac{1}{d}}}{\sqrt{d}-\sqrt{1-\varepsilon_d^2}}\quad.\notag
\end{align}
We apply (\ref{EpsilonFormula}) to \(\rho_d\) and we get

\begin{align}
\frac{\varepsilon_d\sqrt{d}+\sqrt{1-\varepsilon_d^2}\sqrt{1-\frac{1}{d}}}{\sqrt{d}-\sqrt{1-\varepsilon_d^2}}
=\frac{O(\frac{1}{\sqrt{d}})+1+O(\frac{1}{d})}{\sqrt{d}+O(1)}=\frac{1}{\sqrt{d}}+O\left(\frac{1}{d}\right)\quad.\notag
\end{align}
Simplices have only one affine point and hence, the John point and the centroid of \(\tilde{\Delta}_{d+1}\) coincide. It is a classical result for simplices that a line segment from a vertex to the opposite base through the centroid is cut by the centroid in the ratio dimension to \(1\). Therefore, the centroid of \(\tilde{\Delta}_{d+1}\) is 

\begin{align}
\left(\frac{1}{d+2}(-\rho_d)+\frac{d+1}{d+2}\sqrt{1-\frac{1}{d}}\right)e_1=\left(\left(1-\frac{1}{d+2}\right)\sqrt{1-\frac{1}{d}}+O(\frac{1}{d^{3/2}})\right)e_1\quad.\notag
\end{align}
The John ellipsoid is symmetric about the centroid and, at least if \(d\) is large enough, we conclude that 
the John ellipsoid is included in 

\begin{align}
\mathrm{conv}\left[(-\varepsilon_d,\sqrt{1-\varepsilon_d^2} \Delta_d),  (\sqrt{1-\frac{1}{d}}, \sqrt{d}\Delta_d)\right]\quad.\label{TruncSimplex}
\end{align}
Hence, it is also the John ellipsoid of (\ref{TruncSimplex}), by the preceding discussion also of \(K_d\) and \(C_d\). We are ready to calculate \(\phi_{j,l}(C_d)\)

\begin{align}
1-\phi_{j,l}(C_d)
=&\frac{\|j(C_d)-l(C_d)\|_2}{\mathrm{vol}_1(a\cap C_d)}
=\frac{\left(1-\frac{1}{d+2}\right)\sqrt{1-\frac{1}{d}}+O\left(\frac{1}{d^{3/2}}\right)}{\sqrt{1-\frac{1}{d}}+\varepsilon_d}\notag\\
=& \frac{\left(1-\frac{1}{d+2}\right)\left(1-\frac{1}{2d}+O(\frac{1}{d^2})\right)+O\left(\frac{1}{d^{3/2}}\right)}{1-\frac{1}{2d}+O(\frac{1}{d^2})+\frac{1}{d}+O(\frac{1}{d^2})}\notag\\
=&\frac{1-\frac{3}{2d}+O\left(\frac{1}{d^{3/2}}\right)}{1+\frac{1}{2d}+O\left(\frac{1}{d^2}\right)}
=1-\frac{2}{d}+O\left(\frac{1}{d^{3/2}}\right)\notag\\
=&1-\frac{2}{d+1}(1+o(1))\notag
\end{align}
and this shows the right-hand side of inequality (\ref{NewOrderOfMagnitude}).

\hfill \(\Box\)

\bibliography{ijmsample}
\bibliographystyle{ijmart}

\end{document}